\documentclass[10pt]{article}
\usepackage{hyperref}
\usepackage{amssymb}
\usepackage{amsmath}
\input{liemacs10.sty} 
\numberwithin{equation}{section}

\usepackage[usenames,dvipsnames]{xcolor}
\usepackage{cancel}

\usepackage{color}

\renewcommand{\up}{{\mathop{\uparrow}}}
\renewcommand{\down}{{\mathop{\downarrow}}}

\newcommand{\sH}{{\sf H}}
\newcommand{\sN}{{\sf N}}
\newcommand{\sK}{{\sf K}}
\newcommand{\sV}{{\tt V}}

\newcommand{\RR}{\mathbb{R}}

\renewcommand{\phi}{\varphi} 

\newcommand{\AdS}{\mathop{{\rm AdS}}\nolimits}

\newcommand{\Stand}{\mathop{{\rm Stand}}\nolimits}

\newcommand{\Stab}{{\mathrm{Stab}}}

\renewcommand\mlabel{\label}

\begin{document}
\title{A family of non-modular covariant AQFTs}

\author{
{\bf Vincenzo Morinelli}\\
Dipartimento di Matematica, Universit\`a di Roma ``Tor
Vergata''\\ 
E-mail: {\tt morinell@mat.uniroma2.it}
\and
{\bf Karl-Hermann  Neeb}\\
Department  Mathematik, FAU Erlangen-N\"urnberg, \\ 
E-mail: {\tt  neeb@math.fau.de}}

\date{\today} 

\maketitle 
\begin{abstract}
Based on the construction provided in our paper 
``Covariant homogeneous nets of standard subspaces'', 
Comm. Math. Phys. 386 (2021), 305--358, 
we construct {non-}modular covariant one-particle nets on the two-dimensional de
Sitter spacetime and on the three-dimensional  Minkowski space.
\end{abstract}
\tableofcontents 

\section{Introduction}

The modular theory of von Neumann algebras, 
based on the Tomita--Takesaki Theorem 
 has applications in 
many different areas of mathematics and physics because it 
generates rich structures from  very basic assumptions, 
such as the existence of cyclic separating vectors. 
In Algebraic Quantum Field Theory (AQFT), 
the verification of the Bisognano--Wichmann property was 
a fundamental breakthrough \cite{BW76}. 
This property deeply  relates the geometry of the models to the modular structure of von Neumann algebras of observables.

The description of the QFT models through operator algebras is based on 
an assignment of von Neumann algebras $A(\cO)$ to spacetime regions $\cO$ satisfying isotony.
Basic features of QFT can be expressed as natural properties of such a 
Haag--Kastler net $\cO \mapsto A(\cO)\subset\mathcal B(\cH)$ 
(\cite{Ha96}): 
Algebras associated with spacelike separated regions are required to commute (locality), and covariance is formulated in terms of a 
unitary positive energy representation of the symmetry group (the Poincar\'e group on the Minkowski spacetime or the Lorentz group on de Sitter spacetime) that acts covariantly on the observable algebras. Together with an invariant vacuum vector state, 
this defines the vacuum representation of a net of observables as introduced by Haag and Kastler in 1964 (\cite{HK64}). 

Given a von Neumann algebra $\cM$ 
 with a cyclic and separating vector, 
the Tomita--Takesaki Theorem provides a unitary 
one-parameter group,  called modular group, 
whose adjoint action on the algebra defines 
a one-parameter group of automorphism. The modular groups
of some algebras with 
particular localization properties actually correspond to 
global geometric symmetries. This is 
called the Bisognano--Wichmann (BW) property. It marked one of the formidable 
successes of modular theory: The modular structure of some observable algebras 
in the vacuum sector has geometric meaning.

The BW property  has been verified for large number of models (see
e.g.~\cite{DM20, Mu01, BGL93}) and applied in various ways with feedback both for mathematics 
and for physics. For recent developments concerning entropy of QFT's we refer to \cite{LW22, W18, Lo19, Xu20, LMo20, CLRR22} and for some new constructions exploiting geometric symmetries and modular theory to \cite {LMPR19, MR20}. 
The property of so-called ``modular covariance'' of a net is weaker 
than the Bisognano--Wichmann property \cite{Bo92, Yngvason, BDFS00}.  
On one hand, the Bisognano--Wichmann property assumes covariance with respect to a Poincar\'e (on Minkowski spacetime) or Lorentz (on de Sitter spacetime) representation  and states that the unitary representation of the Lorentz boosts coincides with the modular group of the algebra associated to a wedge domain, properly parametrized.
On the other hand, modular covariance
  ensures that the modular groups of wedge region algebras act
  on the net geometrically as the Lorentz boosts, but not that they belong to the given representation of the symmetry group.
  On Minkowski spacetime the modular group of the wedge regions generate  a representation of the Poincar\'e group satisfying the positive energy condition \cite{GL95}.
  This result does not directly extend to de Sitter spacetime because a specific domain for the modular operator of orthogonal wedges has not been proved yet to be dense in the case of complementary series representations of $\PSL_2(\RR)$ as a subgroup of the Lorentz group (cf.\ \cite[\S 4.4]{MN21}).

A long standing question is how to characterize models where just modular covariance fails or, more generally, that are not modular covariant. 
Most of the models violating the Bisognano--Wichmann property are still modular covariant, see \cite{LMR16, Mo18, Do10}. A first family of models without modular covariance in Minkowski space with dimension $1+d>2$ has been provided by Yngvason in \cite{Yngvason}. {For
  a certain two-point function, he constructs translation
  covariant nets of von Neumann algebras with positive energy, computes
  the modular operators of wedge algebras and concludes} that
they do not act covariantly on the net. {However,}
these models are not expected to be Lorentz covariant. 

A new approach to  geometric features in AQFT models has been provided in 
\cite{MN21}, where wedge regions 
are replaced by  abstract data associated 
to a graded Lie group, which in representations correspond to 
modular operators and conjugations as 
produced by the Tomita--Takesaki Theorem. 
This Lie theoretic perspective on the wedge--boost 
correspondence creates a means to construct  
AQFT models. In \cite{MN21} we also determine a class of Lie groups 
which are compatible with a notion of wedge localization, 
including the known cases from physics 
(see also \cite{NO21, NO22a, NO22b} for related recent work in 
this direction).

The construction of the models is based on the Brunetti--Guido--Longo 
(BGL) construction of the free fields: Nets of von Neumann algebras are constructed through the second quantization canonical procedure starting 
from one-particle nets of standard subspaces defined 
by (anti-)unitary representations of graded Lie groups 
(\cite{BGL02}). 
Thus, in this context, in order to {obtain}
a net of von Neumann algebras,  
one has to construct 
nets of standard subspaces on the one-particle Hilbert space.

In this paper we present a {structural} condition
{that can be used to} construct non-modular covariant nets
in the abstract setting of \cite{MN21}.
Then we construct a family of Lorentz and Poincar\'e covariant models 
based on one-particle nets that are non-modular covariant, 
on  two-dimensional 
de Sitter spacetime and on three-dimensional Minkowski spacetime. 
 The idea is to construct a BGL one-particle net of real subspaces 
from a representation of a ``large'' group $G$ that has a restriction 
to the de Sitter or to the Minkowski symmetry group which is not modular covariant.
It is central that the Lie algebra $\fg$ of 
$G$
contains non-symmetric Euler elements. Several groups 
of this kind are described in \cite{MN21}. 
With the proper identifications, we first obtain non-local nets on
{two-dimensional} de 
Sitter spacetime which are Lorentz covariant. Then
we  generalize the construction to three-dimensional Minkowski space.

The paper is organized as follows.
In Section~\ref{sec:2}  we present the construction of the generalized AQFT models provided in \cite{MN21}.  In Section~\ref{sec:3} 
we present a general construction of non-modular covariant nets 
and give sufficient criteria for its applicability. In Subsection~\ref{ex:SLn} we present an explicit example of this method. 
In Subsection~\ref{sect:minkcount} we show how the construction can be applied to 
Minkowski space. {An outlook on disintegration of covariant representations, locality and higher dimensional examples is given in Section~4.}

\section{Preliminaries}
\mlabel{sec:2}

In this section we collect background on several 
concepts and their properties: 
standard subspaces, abstract Euler wedges and the BGL construction. 
In Subsection~\ref{subsec:2.4} we prepare the group theoretic background 
for our construction of non-modular covariant nets. 

\subsection{One-particle subspaces}

We call a closed real subspace $\sH$ of the complex Hilbert space 
$\cH$ \textit{cyclic} if $\sH+i\sH$ is dense in $\cH$, \textit{separating} if $\sH\cap i\sH=\{0\}$, and \textit{standard} 
if it is cyclic and separating. The symplectic ``complement'' 
of a real subspace $\sH$ is defined by the symplectic form 
$\Im \langle\cdot,\cdot\rangle$ via 
\[ \sH'=\{\xi\in\cH:(\forall \eta \in \sH)\ 
 \Im\langle\xi,\eta \rangle=0 \}.\] 
Then $\sH$ is separating if and only if $\sH'$ is cyclic, hence $\sH$ is standard if and only if $\sH'$ is standard.
For a standard subspace $\sH$, we 
define the {\it Tomita operator} as the closed antilinear involution
\[S_\sH:\sH+i\sH \to \sH + i \sH,\quad 
\xi+i\eta\mapsto \xi-i\eta. \] 
The polar decomposition $S_\sH=J_\sH\Delta_\sH^{\frac12}$ defines an 
antiunitary involution $J_\sH$ (a conjugation) 
and the modular operator~$\Delta_\sH$. 
For the modular group  $(\Delta_\sH^{it})_{t \in \R}$, 
we then have 
\[  J_\sH\sH=\sH', \quad 
 \Delta^{it}_\sH\sH=\sH \quad{\mbox{ and }\quad J_H\Delta^{it}_HJ_H=\Delta^{it}_H}\qquad \mbox{ for every  } \quad 
t\in \R.\]
One also has $\sH=\ker(S_\sH-\textbf{1})$ (\cite[Thm.~3.4]{Lo08}). This construction 
leads to a one-to-one correspondence between Tomita operators and 
standard subspaces: 

\begin{proposition}\label{prop:11}{\rm (\cite[Prop.~3.2]{Lo08})} 
The map $\sH\mapsto S_\sH$ 
is a bijection between the set of standard subspaces of $\cH$ 
and the set of closed, densely defined, antilinear involutions on $\cH$. 
Moreover, polar decomposition $S = J\Delta^{1/2}$ 
defines a one-to-one correspondence between such involutions and 
pairs $(\Delta, J)$, where $J$ is a conjugation and $\Delta >0$ 
selfadjoint with $J\Delta J =\Delta^{-1}$. 
\end{proposition}
The modular operators of symplectic complements satisfy the following relations 
\[ S_{\sH'}=S_{\sH}^*, 
\qquad \Delta_{\sH'}=\Delta_\sH^{-1}, \qquad J_{\sH'}=J_{\sH}.\] 

From Proposition~\ref{prop:11} we easily deduce: 
{\begin{lemma}\label{lem:sym}{\rm(\cite[Lemma 2.2]{Mo18})}
Let $\sH\subset\cH$ be a standard subspace  and $U\in\AU(\cH)$ 
be a {unitary or anti-unitary} operator. 
Then $U\sH$ is also standard and 
$U\Delta_\sH U^*=\Delta_{U\sH}^{\eps(U)}$ and $UJ_\sH U^*=J_{U\sH}$, 
where $\eps(U) = 1$ if $U$ is unitary and 
$\eps(U) = -1$ if it is antiunitary. 
\end{lemma}}

\subsection{Euler wedges}
Let $G$ be  a finite dimensional  $\Z_2$-{\it graded Lie group} 
$(G,\eps_G)$, i.e., $G$ is a Lie group and $\eps_G \:  G \to \{\pm 1\}$ 
a continuous homomorphism. We  write 
\[ G^\up = \eps_G^{-1}(1) \quad \mbox{ and } \quad G^\down = \eps_G^{-1}(-1),\] 
so that $G^\up \trile G$ is a 
normal subgroup of index $2$ and $G^\down = G \setminus G^\up$. 
As the subgroup $G^\up$ is open and closed, 
it contains the connected component $G_e$ of the 
neutral element~$e$ in~$G$. 

\begin{defn} \mlabel{def:euler}
 (a) We call an element $x$ of the finite dimensional 
real Lie algebra $\g$ an
{\it Euler element} if $\ad x$ is non-zero and diagonalizable with 
$\Spec(\ad x) \subeq \{-1,0,1\}$. In particular the eigenspace 
decomposition with respect to $\ad x$ defines a $3$-grading 
of~$\g$: 
\[ \g = \g_1(x) \oplus \g_0(x) \oplus \g_{-1}(x), \quad \mbox{ where } \quad 
\g_\nu(x) = \ker(\ad x - \nu \id_\g)\] 
Then $\sigma_x(y_j) = (-1)^j y_j$ for $y_j \in \g_j(x)$ 
defines an involutive automorphism of $\g$. 

We write $\cE(\g)$ for the set of Euler elements in~$\g$. 
The orbit of an Euler element  $x$ under the group 
$\Inn(\g) = \la e^{\ad \g} \ra$ of 
{\it inner automorphisms} 
is denoted with $\cO_x = \Inn(\g)x \subeq \g$.\begin{footnote}
{For a Lie subalgebra $\fs \subeq \g$, we write 
$\Inn_\g(\fs)= \la e^{\ad \fs} \ra \subeq \Aut(\g)$ for the subgroup 
generated by $e^{\ad \fs}$.}  
\end{footnote}
We say that $x$ is {\it symmetric} if $-x \in \cO_x$.

\nin (b) {The set 
\[ \cG {:= \cG(G)} := \{ (x,\sigma)\in\g \times G^\down \: \sigma^2 = e, \Ad(\sigma)x = x\}\] 
is called the {\it abstract wedge space of $G$}. 
We assign to 
$W = (x,\sigma) \in \cG$ the one-parameter  group 
\begin{equation}
  \label{eq:deflambdaw}
 \lambda_W \: \R \to G^\up \quad \mbox{ by } \quad 
\lambda_W(t) := \exp(t x).
 \end{equation} }
An element $(x,\sigma) \in \cG$ is called an 
{\it Euler couple} or {\it Euler wedge}   
if 
\begin{equation}\label{eq:eul}
\Ad(\sigma)=\sigma_x := e^{\pi i \ad x}.\end{equation} 
Then $\sigma$ is called an {\it Euler involution}. 
We write $\cG_E\subeq \cG$ for the subset of Euler couples.

\nin (c) For the graded group $(G,\eps)$, we consider on $\g$ the 
{\it twisted adjoint action} which changes the sign on odd group elements: 
\begin{equation}
  \label{eq:adeps}
  \Ad^\eps \: G \to \Aut(\g), \qquad 
\Ad^\eps(g) := \eps_G(g) \Ad(g).
\end{equation}
It extends to an action of $G$ on $\cG$  by 
\begin{equation}
  \label{eq:cG-act}
 g.(x,\sigma) := (\Ad^\eps(g)x, g\sigma g^{-1}).
\end{equation}

\nin (d) (Duality operation) 
We define the notion of a  
``causal complement'' on the abstract wedge space: 
For $W = (x,\sigma) \in \cG$, we define the {\it dual wedge} by  
$W' := (-x,\sigma) {= \sigma.W}$. 
Note that $(W')' = W$ and $(gW)' = gW'$ for $g \in G$ 
by \eqref{eq:cG-act}.  

{The relation $\sigma.W = W'$ is our main 
motivation to work with the twisted adjoint action. 
This relation fits the geometric interpretation in the context
of wedge domains in spacetime manifolds.} 

\nin (e) (Order structure on $\cG$) For a given invariant closed 
convex cone $C \subeq \g$, we obtain an order 
structure on $\cG$ as follows. We associate to $W = (x,\sigma) \in \cG$ 
a semigroup $\cS_W$ whose unit group is 
$\cS_W \cap \cS_W^{-1} = G^\up_{W}$, the stabilizer of $W$
 (\cite[Thm.~III.4]{Ne23}). 
It is specified by 
\[ \cS_W := \exp(C_+) G^\up_{W} \exp(C_-) 
= G^\up_{W} \exp\big(C_+ + C_-\big),\] 
where the convex cones $C_\pm$ are the following intersections 
\[ C_\pm := \pm C \cap \g^{-\sigma} \cap \ker(\ad x \mp \1)
\quad \mbox{ and } \quad \g^{\pm \sigma} 
:= \{ y \in \g \: \Ad(\sigma)(y) = \pm y \}. \]
Then 
$\cS_W$ defines a $G^\up$-invariant partial order on the orbit 
$G^\up.W \subeq \cG$ by 
\begin{equation}
  \label{eq:cG-ord}
g_1.W \leq g_2.W  \quad :\Longleftrightarrow \quad 
g_2^{-1}g_1 \in \cS_W.
\end{equation}
In particular, $g.W \leq W$ is equivalent to $g \in \cS_W$. 

\end{defn} 

\begin{lem} \mlabel{lem:1.1} {\rm(\cite[Lemma~2.6]{MN21})}
For every $W=(x_W,\sigma_W) \in \cG$, $g \in G$, and $t \in \R$, 
the following assertions hold: 
  \begin{itemize}
\item[\rm(i)] $\lambda_W(t).W = W, \lambda_W(t).W' = W'$ and $\sigma_W.W = W'.$ 
\item[\rm(ii)] $\sigma_{W'} = \sigma_W$ and 
$\lambda_{W'}(t) = \lambda_W(-t)$. 
\item[\rm(iii)] $\sigma_W$ commutes with $\lambda_W(\R)$. 
\end{itemize}
\end{lem}

\begin{rem} \mlabel{rem:stab} Let $W = (x,\sigma) \in \cG$ 
and consider $y \in \g$. Then $\exp(\R y)$ fixes $W$ if and only if 
\[ [y,x] = 0 \quad \mbox{ and  } \quad y = \Ad(\sigma)y.\] 
If $(x,\sigma)$ is an Euler couple, then 
$\Ad(\sigma)y = e^{\pi i \ad x}y =  y$ follows from $[y,x] =0$, so that 
\begin{equation}
  \label{eq:centralizercond}
\g_W := \{ y \in \g \:  \exp(\R y) \subeq G_W^\up\}  
= \ker(\ad x).
\end{equation}
\end{rem}

\begin{defn} {\rm(The abstract wedge space)}
  \mlabel{def:2.6}
From here on, we always assume that $\cG \not=\eset$, i.e., that 
$G^\down$ contains an involution $\sigma$. Then 
\[ G \cong G^\up \rtimes \{\id,\sigma\} \]
For a fixed couple $W_0 = (h,\sigma) \in \cG$, the orbits 
\[ \cW_+ (W_0):= G^\up.W_0 \subeq \cG  \quad \mbox{ and } \quad 
\cW(W_0) := G.W_0 \subeq \cG  \] 
are called the {\it positive} and the {\it full wedge space containing $W_0$}. 
\end{defn}

\begin{rem} \label{ex:desit}
(Lorentz wedges on de Sitter spacetime) 
The de Sitter spacetime is the manifold 
$\dS^{d}=\{(t,\bx)\in\RR^{1+d}: \bx^2-t^2=1\}$, 
endowed with the metric obtained by restriction of 
the Minkowski metric 
\[ ds^2=dt^2-dx_1^2-\ldots-dx_d^2 \] 
 to $\dS^d$. 

The generator $k_{1} \in \so_{1,d}(\R)$ of the Lorentz boost on the 
$(x_0,x_1)$-plane 
\[  k_1(x_0,x_1,x_2, \ldots, x_{d}) = (x_1, x_0, 0, \ldots, 0)\] 
is an Euler element. It combines with the spacetime 
reflection 
\[ j_1(x) =(-x_0,-x_1,x_2,\ldots,x_d) \] 
to the Euler couple $(k_1, j_1)$ for the graded 
Lie group $\SO_{1,d}(\R)$. 
The spacetime region 
\[ W_{x_1}=\{x\in\RR^{1+d}: |x_0|<x_1\}, \] 
is called the {\it standard right wedge} and we put 
\[ W_{x_1}^{\dS} := W_1 \cap \dS^d. \] 
Note that $W_{x_1}$ and therefore $W_{x_1}^{\dS}$ are 
invariant under $\exp(\R k_{1})$. 
Lorentz transforms {$W^{\dS}=g.W_1^{\dS}$} of $W_{x_1}^{\dS}$  with $g\in \SO_{1,d}(\R)$ 
are called {\it wedge regions in de Sitter space}. 
They are in 1-1 correspondence with  
Euler couples in $\cG(\SO_{1,d}(\R)$ 
 and one can associate to $W_{x_1}$ 
the couple $(k_1, j_1) \in \cG_E(\SO_{1,d}(\R))$ (cf.~\cite[Lemma~4.13]{NO17} 
and \cite[Sect5.2]{BGL02}).
{For $\ell = 2,\ldots, d$, one likewise obtains couples
$(k_\ell, j_\ell)$ corresponding to the wedges $W_{x_\ell}
= \{ x \in \R^{1+d} \: |x_0| < x_\ell\}$.}

In this paper we will focus on 2-dimensional de Sitter spacetime. {The orbit of the wedge $W_{x_1}^{\dS}$ under the Lorentz group is indicated with $\cW^{\dS}_+=\cL_+^\up W_{x_1}$ where $\cL_+^\up:=\SO_{1,2}(\R)^\up$}. In this case the orthochronous symmetry group $\SO_{1,2}(\R)^\up$ 
is isomorphic to $\PSL_2(\RR)=\SL_2(\RR)\slash\{\pm1\}$
and the action on $\dS^2$ can be visualized by the adjoint action of 
$\PSL_2(\RR)$ on the orbit generated by the matrix $\frac12 \diag (1,-1)$. 
{We shall use the following matrix picture of $\dS^2$,
  which is implemented by the bijection}
\begin{align*} \label{eq:desitterid}
\dS^2 &\rightarrow{\dS^2_{\rm mat}}:=\Big\{ X\in M_2(\RR) \: \tr X = 0, \det X =-\frac14\Big\} \subeq \fsl_2(\R)\\
 x=(x_0,x_1,x_2) &\mapsto \tilde x:=\frac12\left(\begin{array}{cc}x_1&-x_0-x_2\\x_0-x_2&-x_1\end{array}\right), 
  \end{align*} %x_0\sigma_0+x_1\sigma_2-x_2\sigma_1

  and
\[ \sigma_0=\frac12\left(\begin{array}{cc}0&-1\\1&0
\end{array}\right),\; 
\sigma_1=\frac12\left(\begin{array}{cc}0&1\\1&0
\end{array}\right),\;
\sigma_2=\frac12\left(\begin{array}{cc}1&0\\0&-1
\end{array}\right).
\]
Note that, for
$X \in\dS^2_{\rm mat}$, 
the vector $y=({-2}\Tr(X\sigma_0),{2}\Tr(X\sigma_2),-
{2}\Tr(X\sigma_1))\in \dS^2$ satisfies 
$\tilde y = X$.
{The Minkowski quadratic form $x^2 = x_0^2 - x_1^2 - x_2^2$
  corresponds to the determinant by}  
$x^2=4\det \tilde x$, so that 
$x\in \dS^2$ if and only if $\det\tilde x=-\frac14$. 

We write $\Lambda: \SL_2(\RR)\rightarrow \cL_+^\up$ for the quotient map
defined by the relation: 
\begin{equation}\label{eq:Lambdaact}
(\Lambda(g)x)\,\tilde{}=g\tilde x g^{-1} \quad \mbox{ for } \quad x\in \dS^2, 
g\in \SL_2(\RR).
\end{equation} 
Then it is easy to see that $\dS^2_{\rm mat}=\{g\sigma_1g^{-1}:g\in\SL_2(\RR)\}$. 
The one-parameter groups 
\begin{equation}
  \label{eq:lambda-xi}
 \lambda_{x_i}(t)=\exp(\sigma_it)\in\SL_2(\RR), i =1,2,  
\end{equation}
are the lifts of the boosts
$\Lambda_{W_{x_i}}(t)\in\cL_+^\up$, and $r(\theta)=\exp(-\sigma_0\theta)$ is the one-parameter group lifting the 
one-parameter group 
\begin{equation}
  \label{eq:R}
 \Lambda(r(\theta)) = R(\theta)=\pmat{ 1 & 0 & 0 \\ 
0 & \cos \theta&-\sin\theta\\ 
0 & \sin\theta&\cos\theta}
\end{equation}
of space rotations.
 \end{rem}

\begin{theorem} \mlabel{thm:classif-symeuler} {\rm(\cite[Thm.~3.10]{MN21})}
Suppose that $\g$ is a non-compact simple 
real Lie algebra and that  $\fa \subeq \g$ is maximal $\ad$-diagonalizable 
with restricted root system 
$\Sigma = \Sigma(\g,\fa) \subeq \fa^*$ of type $X_n$. 
We follow the conventions of the tables in {\rm\cite{Bo90}}
for the classification of irreducible root systems and the enumeration 
of the simple roots $\alpha_1, \ldots, \alpha_n$. 
For each $j \in \{1,\ldots, n\}$, we consider the uniquely determined element 
$h_j \in \fa$ satisfying $\alpha_k(h_j) =\delta_{jk}$.
Then every Euler element in $\g$ is conjugate under inner automorphism 
to exactly one~$h_j$.For every irreducible root system, 
the Euler elements among the $h_j$ are  the following: 
\begin{align} %\label{eq:eulelts}
&A_n: h_1, \ldots, h_n, \quad 
\ \ B_n: h_1, \quad 
\ \ C_n: h_n, \quad \ \ \ D_n: h_1, h_{n-1}, h_n, \quad 
E_6: h_1, h_6, \quad 
E_7: h_7.\label{eq:eulelts2}
\end{align}
For the root systems $BC_n$, $E_8$, $F_4$ and $G_2$ no Euler element exists 
(they have no $3$-grading). 
The symmetric Euler elements are 
\begin{equation}
  \label{eq:symmeuler}
A_{2n-1}: h_n, \qquad 
B_n: h_1, \qquad C_n: h_n, \qquad 
D_n: h_1, \qquad 
D_{2n}: h_{2n-1},h_{2n}, \qquad 
E_7: h_7.  
\end{equation}
\end{theorem}

{\begin{rem}  \label{rem:2.9}
The preceding theorem shows that non-symmetric 
Euler elements exist 
for the root systems of type $A_n, n \geq 2$, $D_n, n \geq 4$, 
and $E_6$.  
\end{rem}}

\begin{ex} \mlabel{ex:2.9} For $\g = \fsl_{n}(\R)$, 
the subspace 
\[ \fa = \Big\{ \diag(x_1,\ldots, x_n) \:  \sum_j x_j = 0\Big\}\]
of diagonal matrices is maximal abelian. 
In terms of the linear functionals 
$\eps_j(\diag(x)) = x_j$ on $\fa$, the root system is 
\[ A_{n-1} = \{ \eps_i - \eps_j \: i\not=j \in \{1,\ldots, n\} \}. \] 
The matrices 
\[ h_j = \frac{1}{n} \pmat{ (n-j) \1_j & 0 \\ 0 & - j \1_{n-j}}, 
\quad j = 1,\ldots, 
n-1, \] 
are Euler elements. They are symmetric if and only if~$n = 2j$. 
A corresponding graded Lie group is $G = \PGL_n(\R)$ with 
$G^\up = \PSL_n(\R)$. 
\end{ex}

\subsection{Nets of standard subspaces}
Hereafter we will consider orbits of Euler elements $\cW_+\subset\cG_E(G)$
(Definition~\ref{def:2.6}).

\begin{defn}\label{def:net} 
Let {$G = G^\up \rtimes \{e,\sigma\}$ be as above, 
$C \subeq \g$ be a closed convex $\Ad^\eps(G)$-invariant cone in $\g$, 
and fix a $G^\up$-orbit $\cW_+ \subeq \cG_E(G)$.}
Let $(U,\cH)$ be a unitary representation of $G^\up$  and 
\begin{equation}
  \label{eq:net}
\sN \: \cW_+ \to \Stand(\cH) 
\end{equation}
be a map, also called a {\it net of standard subspaces}. 
In the following we denote this data as $(\cW_+,U,\sN)$.
We consider the following properties: 
\begin{itemize}
\item[\rm(HK1)] {\bf Isotony:} $\sN(W_1) \subeq \sN(W_2)$ for $W_1 \leq W_2$. 
  \begin{footnote}
{Here we refer to the order structure on $\cW_+$ introduced in 
Definition~\ref{def:euler}(e).}
  \end{footnote}
\item[\rm(HK2)] {\bf Covariance:} $\sN(gW) = U(g)\sN(W)$ for 
$g \in {G}^\up$, $W \in \cW_+$. 
\item[\rm(HK3)] {\bf Spectral condition:} 
$ C \subeq C_U := \{ x \in \g \: -i \partial U(x) \geq 0\}$, {where 
$U(\exp tx) = e^{t\partial U(x)}$ for $t \in \R$. }
We then say that $U$ is {\it $C$-positive}.  
\item[\rm(HK4)]{\bf Locality:} 
If $W \in \cW_{+}$ is such that $W'\in\cW_+$, then $\sN(W') \subset \sN(W)'.$
\item[\rm(HK5)] {\bf Bisognano--Wichmann (BW) property:} 
$U(\lambda_W(t)) = \Delta_{\sN(W)}^{-it/2\pi}$ for all 
$W \in \cW_+, t \in \R$.
\item[(HK6)] \textbf{Haag Duality}: $\sN(W{'})= \sN(W)'$ 
{for all $W \in \cW_+$ with $W' \in \cW_+$}.
\item[(HK7)] \textbf{G-covariance:} 
There exists an (anti-)unitary extension of $U$ from $G^\up$ to $G$ such that 
\begin{equation}
  \label{eq:twitcovar}
 \sN(g. W) = U(g) \sN(W) \quad \mbox{ for } \quad 
g \in G, W \in \cW_+.
 \end{equation}
\item[(HK8)] {\bf PCT property:}  
Suppose that (HK7) is satisfied and that $U$ is the corresponding 
representation. Then 
$U(\sigma_W)=  J_{\sN(W)}$ for  $W \in \cW_+$ {with} $W^{'}\in \cW_+$. 
\end{itemize}
\end{defn}

\begin{thm} {\rm(Brunetti--Guido--Longo (BGL) net {generalization}, \cite{MN21})}
  \label{thm:BGL} 
If $(U,G)$ is an \break (anti-)unitary representation,  
then we obtain a $G$-equivariant map 
$\sN_U \:  \cG \to \Stand(\cH)$ 
determined for $W = (k_W, \sigma_W)$ by 
\begin{equation}
  \label{eq:bgl}
J_{\sN_U(W)} = 
U(\sigma_W) \quad \mbox{ and } \quad \Delta_{\sN_U(W)}^{-it/2\pi} 
= U(\exp t k_W) 
\quad \mbox{ for } \quad t \in \R.
\end{equation}
\end{thm}

The BGL net associates to every wedge $W\in\cG$ a standard subspace $\sN_U(W)$. We shall denote with $(\cW_+,\sN_U,U)$ the restriction of the BGL net to 
the $G^\up$-orbit $\cW_+ \subeq \cG_E(G)$.

\begin{theorem} {\rm(\cite[Thm.~4.12, Prop.~4.16]{MN21})}\label{thm:MN}
The restriction of the BGL net $\sN_U$ 
associated to an (anti-)unitary $C$-positive  
representation $U$ of $G=G^\uparrow\rtimes\{e,\sigma\}$ 
to a $G^\up$-orbit $\cW_+\subeq \cG_E$ satisfies all the axioms 
{\rm(HK1)--(HK8)}. 
\end{theorem}

We are interested in models where the following property fails. 
\begin{itemize}
\item[\rm(MC)] \textbf{Modular covariance}: $\Delta_{\sN(W_a)}^{-it}\sN(W_b)
=\sN(\lambda_{W_a}(2\pi t).W_b)$ for $W_a,W_b\in\cW_+$, $t \in \R$.
\end{itemize}

Modular covariance is an immediate consequence of the Bisognano--Wichmann property. Indeed, {(HK2) and} the  BW property imply 
\[ \Delta_{\sN(W_a)}^{-it}\sN(W_b) 
=U(\lambda_{W_a}(2\pi t))\sN(W_b) =\sN(\lambda_{W_a}({2\pi t}).W_b).\] 

\subsection{Symmetric and non-symmetric Euler elements} 
\mlabel{subsec:2.4}
For the graded group $G := \PGL_2(\R)$ with Lie algebra 
$\g = \fsl_2(\R)$, the abstract wedge space 
$\cW_+$ can be identified with the set $\cE(\g)$ of Euler elements in $\fsl_2(\R)$. 
Since $\ker(\Ad)$ is trivial, for any Euler element $h \in \g$, we have 
\[ G_h = \{ g \in G \: \Ad(g)h = h\} =  G_{(h,\sigma_h)},\] 
and 
\begin{equation}
  \label{eq:gh1}
 \cE(\g) \cong \Ad(G^\up)h \cong G^\up/G^\up_h.
\end{equation}

Depending on the choice of the positive cone in $\fg$ we have wedge spaces with different order structures:
\begin{rem} The (ordered) symmetric space $\cE(\fsl_2(\R))$ can be identified with
  the following spaces:
  \begin{itemize}
  \item For a non-trivial order (corresponding to $C \not=\{0\}$):
    The set of non-dense open intervals in $\bS^1$ 
    (wedge space of the conformal structure on $\bS^1$)
    (\cite[Rem.~2.9(c)]{MN21}) and the set
    of wedge domains in two-dimensional Anti-de Sitter space 
$\AdS^2$ (\cite[\S 11]{NO22b}). 
\item For a trivial order (corresponding to $C =\{0\}$):
  The set of wedge domains in two-dimensional de Sitter space 
  $\dS^2$ (see Remark~\ref{ex:desit}) (\cite[Lemma~4.13]{NO17} and \cite{NO22a}).
  \end{itemize}
\end{rem}

Let $(U,\cH)$ be an (anti-)unitary representation of 
$G$ and $\sV \subeq \cH$ be a standard subspace
with modular objects $(\Delta_\sV, J_\sV)$. 
Then there exists a well-defined $G$-equivariant map 
\[ \cE(\g) \to \Stand(\cH), \quad {\Ad(g)h} \mapsto U(g)\sV \] 
if and only if $G^\up_h$ is contained in the stabilizer group $G_\sV^\up$ of~$\sV$:
\begin{equation}
  \label{eq:2.1}
 G^\up_\sV = \{ g \in G^\up \: U(g)\sV = \sV \} 
= \{ g \in G^\up \: 
U(g) J_\sV U(g)^{-1} = J_\sV, 
U(g) \Delta_\sV U(g)^{-1} = \Delta_\sV \}. 
\end{equation}

For $h := \frac{1}{2}\diag(1,-1) \in \cE(\g)$, the stabilizer 
group in $\PSL_2(\R) = \Ad(G^\up)$ is the adjoint image of 
\[  \SL_2(\R)_h = \Big\{ \pmat{ a & 0 \\ 0 & a^{-1}} \: a \in \R^\times \Big\} 
\cong \R^\times,\] 
hence connected because $-\1 \in Z(\SL_2(\R))$. Therefore 
\[ G^\up_h = \exp(\R h),\] 
and thus $G^\up_h \subeq G_\sV$ is equivalent to 
$U(\exp \R h)$ commuting with $J_\sV$ and $\Delta_\sV$. 

\begin{lem} \label{lem:slgl}
  Let $\g$ be a finite-dimensional Lie algebra and 
$h \in \cE(\g)$ an Euler element. 
If the image of $h$ in the semisimple quotient 
$\g/\rad(\g)$ is non-zero, then there exists a 
Lie subalgebra $\fb \subeq \g$ containing $h$ such that 
\begin{itemize}
\item[\rm(a)] $\fb \cong \fsl_2(\R)$ if and only if $h$ is symmetric, and 
\item[\rm(b)] $\fb \cong \gl_2(\R)$ if $h$ is not symmetric. 
\item[\rm(c)] If $h$ is symmetric, 
then {$\Inn_\g(\fb) \cong \PSL_2(\R)$.}
\item[\rm(d)] If $h$ is not symmetric and $\g$ is simple, then 
$\Inn_\g([\fb,\fb]) \cong \SL_2(\R)$. 
\end{itemize}
\end{lem}
\begin{prf} (a) As all Euler elements in $\fsl_2(\R)$ are symmetric 
($\Inn(\fsl_2(\R))$ acts transitively on $\cE(\fsl_2(\R))$), 
this follows from \cite[Thm.~3.13]{MN21}. 

\nin (b) Suppose that $h$ is not symmetric and 
pick a maximal abelian hyperbolic subspace $\fa \subeq \g$ containing~$h$. 
With \cite[Prop.~I.2]{KN96} we find an $\fa$-invariant 
Levi complement $\fs \subeq \g$. Then 
$\fa_\fs := \fa \cap \fs$ is maximal hyperbolic in $\fs$ 
and $\fa = \fa_\fs + \fz_\fa(\fs)$. We pick a root 
$\alpha \in \Delta(\fs,\fa)$ with $\alpha(h) = 1$ and 
root vectors $x_\alpha \in \fs_\alpha$ and $y_\alpha \in \fs_{-\alpha}$ 
with $h_\alpha := [x_\alpha, y_\alpha] \not=0$.  
We stress that $x_\alpha\in\fs_1(h)$. We use that 
\[ [x_\alpha, y_\alpha] = \kappa(x_\alpha, y_\alpha) a_\alpha, \] 
where $a_\alpha \in \fa$ is the unique element with 
$\alpha(a) = \kappa(a_\alpha,a)$ for all $a \in \fa$,  
and that the Cartan--Killing form
 $\kappa$ induces a dual pairing $\fs_\alpha \times \fs_{-\alpha} \to \R$. 
Then 
\[ \fb_\alpha := \R x_\alpha + \R y_\alpha + \R h_\alpha \cong \fsl_2(\R)\] 
and $[h,\fb_\alpha] \subeq \fb_\alpha$. 
Hence $\fb := \R h +  \fb_\alpha$ is a Lie subalgebra of $\fg$. 
As $h$ is not symmetric, $h \not\in \fb_\alpha$, and therefore 
$\fb \cong \gl_2(\R)$. 

\nin (c) If $h$ is symmetric, then $\fb = [\fb,\fb] \cong \fsl_2(\R)$ by (a)
and the fact that $\fb$ contains an Euler element of $\g$ 
implies that all simple $\fb$-submodules of $\g$ are either trivial 
of isomorphic to the adjoint representation of $\fsl_2(\R)$ 
(consider eigenspaces of $\ad h$). 
This implies that $\Inn_\g(\fb) \cong \PSL_2(\R)$. 

\nin (d) Suppose that $\g$ is simple. 
If $h$ is not symmetric, then the Weyl group reflection 
$s_\alpha$ corresponding to the root $\alpha$ from above satisfies 
\[ s_\alpha(h) 
= h - \alpha(h) \alpha^\vee = h - \alpha^\vee.\] 
As $h$ is not contained in $\R \alpha^\vee \subeq \fb_\alpha$, we have 
$s_\alpha(h) \not\in \R h$. 

The simplicity of $\g$ ensures that the root system $\Delta = \Delta(\g,\fa)$ 
is irreducible and $3$-graded by~$h \in \fa$. 
Therefore 
\[ \Delta_0 := \{\alpha \in \Delta \: \alpha(h)=0\} \] 
spans a hyperplane in $\fa^*$, which coincides with $h^\bot$, and 
thus $\R h = \Delta_0^\bot$ by duality. 
Since $s_\alpha(h)$ is not contained in $\R h$, there exists a 
$\beta \in \Delta_0$ with $\beta(s_\alpha(h)) \not=0$. 
Now $\beta(h) = 0$ implies 
\[ 0 \not= \beta(s_\alpha(h)) = -\beta(\alpha^\vee).\] 
As $s_\alpha(h)$ is an Euler element, we obtain 
$|\beta(\alpha^\vee)| = 1$. Therefore the central element
$e^{\pi i \ad \alpha^\vee}$ of $\Inn_\g(\fb_\alpha)$ acts non-trivially,
and this implies that 
$\Inn_\g(\fb_\alpha) \cong \SL_2(\R)$ because it is a linear 
group with non-trivial center (\cite[Ex.~9.5.18]{HN12}). 
\end{prf}

\begin{remark}\label{rmk:nocent}
Note that the $\fsl_2(\RR)$-subalgebra $\fb_\alpha$ 
generated by $x_\alpha,\,y_\alpha,\,h_\alpha$  does not centralize~$h$. 
We actually have 
\[ h = h_c + \frac{1}{2} \alpha^\vee \quad \mbox{ with } \quad h_c \not=0\ 
\quad \mbox{ and } \quad [h_c,\fb_\alpha] = \{0\}.\]
\end{remark}

\section{{One-particle nets which are not modular covariant}}
\mlabel{sec:3}

Based on the preceding discussion,  
we describe in  {in Subsection~\ref{sect:constr1}}
a general principle that can be used to 
construct one-particle nets that are not modular covariant.
{This is then applied to obtain such nets on two-dimensional
 de Sitter space $\dS^2$ and three-dimensional Minkowski space. }

\subsection{A general  construction principle} \label{sect:constr1}

We describe 
a construction principle for non-modular covariant nets of standard subspaces.
Let $G = G^\up \rtimes \{\1,\tau\}$ be a graded Lie group 
and $(U,\cH)$ an (anti-)unitary representation of $G$. 
We consider the following situation: 
\begin{itemize}
\item a graded subgroup $H \subeq G$, 
\item $W_1 = (h_1, \tau_1)\in \cG_E(H)$ and an Euler couple 
$W_2 = (h_2, \tau_2) \in \cG_E(G)$, so that 
$\Ad(\tau_2) = e^{\pi i \ad h_2}$. 
\item the stabilizer $H^\up_{W_1}$ of $W_1$ in $H^\up$ fixes $W_2$. 
As $\exp(\R h_1) \subeq H^\up_{W_1}$, this implies 
$[h_1, h_2] = 0$ 
(Remark~\ref{rem:stab}). 
\end{itemize}

Then the BGL construction provides a 
standard subspace $\sN_2 = \sN(h_2, \tau_2, U)$ with 
\[ J_{\sN_2} = U(\tau_2) \quad \mbox{ and  }\quad 
\Delta_{\sN_2} = e^{2 \pi i \cdot \partial U(h_2)}.\] 
As the BGL net $\cG(G) \to \Stand(\cH)$ is $G$-equivariant, 
maps $W_2$ to $\sN_2$, and $H^\up$ fixes $W_2$, 
we obtain an $H^\up$-equivariant map 
\[ \sN \: \cG(H) \supeq \cW_+ := 
\cW_+(H,h_1,\tau_1) := H^\up.W_1  \to \Stand(\cH), \quad 
g.W_1 \mapsto U(g)\sN_2\] 
which is uniquely determined by 
\begin{equation}
  \label{eq:incond}
 \sN(W_1) = \sN_2.
\end{equation}

\begin{lem} \mlabel{lem:4.2} 
The net $\sN$ on $\cW_+(H,h_1, \tau_1)$ satisfies  modular covariance 
if and only if, for all ${g \in H^\up},  t \in \R$, the operator 
\[ U(g) U(\exp t (h_1-h_2)) U(g)^{-1} \] 
fixes the standard subspace $\sN_2$, i.e., 
\begin{equation}
  \label{eq:covcond1}
 g \exp(t (h_1-h_2))g^{-1} \in G_{\sN_2} \quad \mbox{ for }\quad 
g \in H^\up, t \in \R. 
\end{equation}
\end{lem}

\begin{prf} The net $\sN$ satisfies  modular covariance if and only if  
\begin{equation}
  \label{eq:2.2}
\sN(\lambda_{g_1 W_1}(t) g_2 W_1) 
\ {\buildrel ! \over =}\ 
\Delta_{\sN(g_1 W_1)}^{-it/2\pi} \sN(g_2 W_1)\quad \mbox{ for } \quad 
g_1, g_2 \in H^\up, t \in \R.
\end{equation}
By covariance of $\sN$, the left hand side equals 
\[ U(\lambda_{g_1 W_1}(t)) \sN(g_2 W_1) 
=  U(g_1 \exp(t h_1) g_1^{-1}) U(g_2) \sN(W_1) 
=  U(g_1) U(\exp(t h_1) U(g_1^{-1} g_2) \sN_2, \] 
and the right hand side is 
\begin{align*}
\Delta_{\sN(g_1 W_1)}^{-it/2\pi} U(g_2) \sN_2
&= \Delta_{U(g_1) \sN_2}^{-it/2\pi} U(g_2) \sN_2 
= U(g_1) \Delta_{\sN_2}^{-it/2\pi} U(g_1)^{-1} U(g_2) \sN_2\\
&= U(g_1) \exp(th_2) U(g_1^{-1} g_2) \sN_2.
\end{align*}
Note that $[h_1, h_2]=0$ implies that 
$U(\exp th_1) U(\exp -th_2) = U(\exp t(h_1-h_2))$.
So \eqref{eq:2.2} means that 
\[   U(g_2^{-1}  g_1) U(\exp -t h_2) U(\exp th_1) U(g_1^{-1} g_2)
=  U(g_2^{-1}  g_1) U(\exp t(h_1-h_2))U(g_1^{-1} g_2)
 \] 
fixes $\sN_2$ for $g_1, g_2 \in H^\up, t \in \R.$  This is \eqref{eq:covcond1}. 
\end{prf}

\begin{rem} \mlabel{rem:4.3}
(a) Condition \eqref{eq:covcond1} is not  easy to evaluate,
but one can easily formulate sufficient conditions for it to be satisfied. 
As $g G_{\sN_2} g^{-1} = G_{U(g)\sN_2}$ by covariance, it is equivalent to 
\begin{equation}
  \label{eq:covcond1b}
 \exp(t(h_1-h_2))\in G_{g.\sN_2} \quad \mbox{ for all }\quad 
g \in H^\up, t \in \R 
\end{equation}
or 
\begin{equation}
  \label{eq:covcond1c}
U(\exp(\R(h_1-h_2))) \subeq 
 \{ U \in \U(\cH) \: (\forall \sH \in \sN(\cW_+)) \ U\sH = \sH\}.
\end{equation}

\nin (b) If $\ker(U)$ is discrete, then the representation of the Lie algebra 
$\g$ is faithful, so that 
any element in the Lie algebra $\g_{\sN_2}$ of the stabilizer group $G_{\sN_2}$ 
commutes with $h_2$, hence is contained in $\g_{W_2} = \ker(\ad h_2)$ 
(cf. Remark~\ref{rem:stab}). We thus have 
\[ \g_{\sN_2} = \g_{W_2} = \ker(\ad h_2).\]
Therefore 
\eqref{eq:covcond1} is equivalent to 
\begin{equation}
  \label{eq:neccond1}
\Ad(H^\up)(h_1 - h_2) \subeq \g_{\sN_2} 
= \ker(\ad h_2),
\end{equation}
which implies in particular, by derivation and by the closedness of 
the Lie subalgebra, that 
\begin{equation}
  \label{eq:incl}
[\fh, h_1 -h_2] \subeq \ker(\ad h_2).  
\end{equation}
We conclude that \eqref{eq:neccond1} is violated if 
\begin{equation}
  \label{eq:negcond}
 [\fh, h_1 - h_2] \not\subeq \ker(\ad h_2).
\end{equation}

\nin (c) If, in addition, $\ad(h_2)\res_{\fh} = - \ad(h_1)$, then  
$[\fh, h_1 - h_2] = [\fh,h_1] \subeq \fh$, so that 
\eqref{eq:negcond} is equivalent to 
$[h_1,\fh] \subeq \ker(\ad h_1)$, which, by 
semisimplicity of $\ad h_1$ on $\fh$, is equivalent to 
$h_1 \in \fz(\fh)$. Therefore $h_1 \not\in \fz(\fh)$ implies \eqref{eq:negcond}. 

\nin (d) The condition under (c) is always satisfied if 
\[ \fh = \g \quad \mbox{ and } \quad h_1 = - h_2 \in \cE(\g).\] 
\end{rem}

\begin{ex}\label{ex:HK} {An easy }example satisfying the construction is the following. Let $U$ be an \break (anti-)unitary representation of the proper 
Lorentz group $\cL_+ = \cL_+^\up\cup \cL_+^\down$ on a Hilbert space 
$\cH_U$ and let $\sN_U$ be the corresponding 
BGL net of standard subspaces.
Fix a wedge $W \in \cG_E(\cL_+)$ 
and consider the net of standard subspaces obtained by 
$\sN_U'(W):=\sN_U(W)'$. Then $\sN_U'$ is a local 
net in the sense that it satisfies (HK4). 
But it does not satisfy modular covariance, since 
$\Delta_{\sN_U'(W)}^{it}=U(\Lambda_W(2\pi t))$ acts with the inverted 
boost flow with respect to the Lorentz action on Minkowski spacetime. In particular it is easy to see that, with the previous notation $h_1=-h_2$.
\end{ex}

\subsection{Non-modular covariant nets}
\subsubsection{On de Sitter space $\dS^2$}\label{sect:desittcex}

Here we apply the previous prescription in order to construct 
nets of standard subspaces {\bf without } 
the modular covariance property.
We can specify the assumption (HK2) for the (double covering) of the Lorentz group as follows. We say that a net of standard subspaces on the
  two-dimensional de Sitter spacetime $W^{\dS}\mapsto \sN(W^{\dS})$ is
  \textbf{Lorentz covariant} with respect to
  a unitary representation $U$ of $\SL_2(\RR)\simeq\widetilde{\mathcal L}_+^\up$
 if \begin{equation}\label{eq:Lorcov}
U(g)H(W^{\dS})= H(\Lambda(g)W^{\dS}) \quad \mbox{ for }\quad g\in\SL_2(\RR), 
\end{equation}
where $\Lambda \: \SL_2(\R) \to \SO_{1,2}(\R)^\up$ is the covering homomorphism, {see \eqref{eq:Lambdaact}.}
We shall say that a net of standard subspaces on the three-dimensional Minkowski spacetime $\RR^{1+2}\supset W^{\RR^{1+2}}\mapsto \sN(W^{\RR^{1+2}})$ is \textbf{Poincar\'e covariant }with respect to a unitary representation $U$  of the (double covering of the) Poincar\'e group $\widetilde{\mathcal{P}}_+^\up\simeq\RR^{1+2}\rtimes_\Lambda\SL_2(\RR)$ if \begin{equation}\label{eq:Poicov}
  U(g)H(W^{\dS})= H(\Lambda(s)W^{\dS}+a)
  \quad \mbox{ for }  \quad g=(a,s)\in\RR^{1+2}\rtimes\SL_2(\RR). 
\end{equation}

Let $G^\up=\Inn(\fg)$ with $\fg$ a simple non-compact Lie algebra  
 as in Theorem~\ref{thm:classif-symeuler} and let $G$ be the graded extension 
$G^\up\rtimes\{1,\sigma\}$, where $\sigma$ is an Euler involution.  We consider 
an {(anti-)unitary} representation 
$U$ of $G$ on a Hilbert space $\cH$.
We pick a {\bf non-symmetric} Euler element $h\in\fg$ 
and the associated wedge $W_h=(h,\sigma_h) \in \cG_E(G)$.
By Lemma \ref{lem:slgl} there exists a 
$\fgl_2$-subalgebra $\fb\subeq \g$ containing $h$  
such that  $\fh := [\fb,\fb] \cong \fsl_2(\R)$ 
satisfies  $[\fh,h]\neq 0$. We consider 
$H^\up := \Inn_\g(\fh) \cong \SL_2(\RR)$ and $H=H^\up\rtimes\{1,\sigma_h\}$
(Lemma~\ref{lem:slgl}(d)). 
{The Euler element $h \in \fb$ has a central component~$h_c$, so that
  \begin{equation}
    \label{eq:hsplit1}
    h = h_c - h_1 \quad \mbox{ with } \quad
h_c \in \fz(\fb) \quad \mbox{ and } \quad
h_1 \in {\cE(\fh).}
  \end{equation}
Choosing the isomorphims $\fh\to\fsl_2(\RR)$ suitably, 
 we may henceforth assume that
\begin{equation}
  \label{eq:hsplit2}
  h_1 = \sigma_1= \frac12\left(\begin{array}{cc}0&1\\1&0
\end{array}\right).
\end{equation}
}

Let $V$ be the restriction $U|_H$. 
The group  $H^\up\cong \SL_2(\RR)$ is the double covering of $\PSL_2(\R) 
\cong \cL_+^\up$. So it 
acts on de Sitter space $\dS^2$ 
through the covering map 
$\Lambda: \SL_2(\RR)\ni g\mapsto \Lambda(g)\in \PSL_2(\RR)$
from Remark~\ref{ex:desit}.

We thus obtain a net $\sH^{\dS}$
on de Sitter spacetime. 
  \begin{theorem} \mlabel{thm:3.5}
  There exists a 
  net of standard subspaces on de Sitter spacetime $\dS^2$ 
     \[ {\cW^{\rm dS}}\ni W^{\rm dS}\mapsto \sH^{\dS}(W^{\rm dS})\subset\cH_U \]
such that
\begin{equation}
  \label{eq:exgen1}
  \sH^{\dS}(\Lambda(g)W_{x_1}^{\rm dS}):=V(g)\sN_U(W_h)\quad 
  \mbox{ for  } \quad g \in H^\up \cong \SL_2(\R),
\end{equation}
where $\sN_U$ is the BGL net defined by $U$, satisfying the Lorentz covariance property \eqref{eq:Lorcov} w.r.t. the representation $V$ of  $H^\up$ and the action defined by \eqref{eq:Lambdaact}.
The net $\sH^{\dS}$ is Lorentz covariant and does not satisfy modular covariance. 
\end{theorem}
We remark that since all the wedge inclusions are trivial and since the
positive cone $C$ is $\{0\}$ for the Lorentz group,
(HK1) and  (HK3) are also satisfied by $\sH^{\dS}$ for trivial reasons. 

\begin{proof} For $\sH^{\dS}$ to be well-defined, we have to argue that
  the stabilizer of $W^{\dS}_{x_1}$ in $\SL_2(\R)$,
  which is the centralizer $H_{h_1}^\up = \{\pm \1\} \exp(\R h_1)$,
  fixes the standard subspace $\sN_U(W_h)$.
  This follows from the fact that $h = h_c - h_1$, where
  $h_c$ is fixed by $H^\up$. Hence $h$ is fixed by $H_{h_1}^\up$, and therefore
  $H_{h_1}^\up$ leaves $\sN_U(W_h)$ invariant by covariance of the
  BGL net $\sN_U$ under $V = U\res_H$.
 
  The covariance of the net $\sH^{\dS}$ follows 
  immediately from its definition in \eqref{eq:exgen1}.
  To see that it is not modular covariant, recall that  the Euler element $h$ is not symmetric, so that it is not contained in any $\fsl_2$-subalgebra 
(\cite[Thm.~3.13]{MN21}),  
and the same holds for all elements in $\Ad(H^\up)h$. 

We put $h_2 := h$ and, as above, we  write $h = h_c - h_1$, where 
$h_c \in \fz(\fb)$ and $h_1 \in \fh = [\fb,\fb]$
is an Euler element in $\fh$. 
Then 
$h_2 - h_1 = h_c - 2 h_1$ satisfies
\[ [h_2-h_1,\fh] = [h_1, \fh] \not\subeq \ker(\ad h_1),\]
so that Remark~\ref{rem:4.3}(b),(c) imply  that the net is not 
modular covariant. More concretely, 
\[ \Ad(H^\up)(h_2 -h_1)
=  h_c - 2 \Ad(H^\up) h_1 \in h_c - 2 \cE(\fh) \]
is not contained in the centralizer of $h = h_2$
(cf.\ \eqref{eq:neccond1}). So the action of the modular group
$U(\exp(-2\pi th)) = \Delta_{\sH^{\dS}(W_{x_1}^{\dS})}^{it}$
on the net differs from that of the one-parameter group
$U(\exp(-2\pi t h_1))$ of the corresponding boost.
\end{proof}

\subsubsection{An explicit example on de Sitter space $\dS^2$}
\label{ex:SLn}
In this section we will present an explicit construction in $\PSL_n(\RR)$ for the previous construction. 

With the previous notation we specify the case 
$\PGL_n(\R)$ as an example. 
Let $G^\up := \PSL_n(\R)$ with Lie algebra
$\g = \fsl_n(\R)$, so that $G^\up \cong \Inn(\g)$. The conjugacy classes of Euler elements in 
$\g$ are represented by the diagonal matrices 
$h_j$, $j = 1,\ldots, n-1$, from Example~\ref{ex:2.9}. 
As the dimensions of the eigenspaces 
show, no two of these elements are conjugate.
This follows from the classification for the root system $A_{n-1}$ 
(Theorem~\ref{thm:classif-symeuler}). 
The Euler element $h_j$ is symmetric only if $n$ is even and 
$j = \frac{n}{2}$. 

More generally, there exist non-symmetric Euler elements
which are diagonal matrices of the form 
\begin{equation}
  \label{eq:dag}
 h=\diag (a_1,\ldots, a_n)\quad \mbox{  with } \quad 
\sum_{i=1}^n a_i=0 \quad \mbox{ and } \quad 
a_i - a_j \in \{-1,0,1\} \quad \mbox{ for } \quad i \not=j.
\end{equation}
We fix such a non-symmetric Euler element $h$. 
After a permutation of the entries, we may 
assume  that {$a_2 - a_1 = 1$, 
  so that $h$ does not commute with the
  subalgebra $\fh \cong \fsl_2(\R)$, generated by the Euler elements 
\[ k_1 = \frac{1}{2} \pmat{ 
1 & 0 &  0\\ 
0 & -1 & 0 \\ 
0 & 0 & \0_{n-2}} \quad \mbox{ and }\quad 
k_2 = \frac{1}{2} \pmat{ 
0 & 1 & 0 \\ 
1 &  0 & 0 \\ 
0 & 0 & \0_{n-2}}. 
\]
Then
\begin{equation}
  \label{eq:hdec3}
  h = h_c - k_1 \quad \mbox{ with } \quad 0\not= h_c
  = \diag(a_1+{\textstyle{\frac{1}{2}}},a_1+{\textstyle{\frac{1}{2}}}
  , a_3, \ldots, a_n).
\end{equation}

Let $\sigma_h \in \Aut(\g)$ be the involution defined by $h$ and
$G := G^\up \{\1,\sigma_h\} \subeq \Aut(\g)$. 
Let $H \subset G$ be the graded subgroup  generated by
the two one-parameter groups $\exp(\R k_{1,2})$ and~$\sigma_h$.
The Lie subalgebra $\fb \subeq \g$ generated by 
$h$, $k_1$ and $k_2$,} is easily seen to be isomorphic to 
$\gl_2(\R)$ with commutator algebra $\fh = [\fb,\fb] \cong \fsl_2(\R)$.
By Lemma~\ref{lem:slgl}(d), we have $H^\up \cong \SL_2(\RR)$,
realized via
\[ \SL_2(\R)\ni g\mapsto
g_{\PSL_n}:=\Big[\pmat{g&0\\0&\1_{n-2}}\Big]\in \PSL_n(\R).\]
Observe that $\sigma_h\res_\fh = \sigma_{k_1}$. 
Let $W_h=(h,\sigma_h) \in \cG_E(G)$ be the Euler wedge
associated to~$h$ and observe that its stabilizer in
$H^\up$ coincides with the stabilizer of $(k_1, \sigma_h) \in \cG(H)$.

Let $(U,\cH_U)$ be an (anti-)unitary representation 
of $G = \PGL_n(\R)$, so that $U$ is unitary on $G^\up=\PSL_n(\RR)$
(see \cite[Lemma~2.10]{NO17} for existence).  
In order to fit with Remark~\ref{ex:desit} (note that $k_1$ corresponds
to $\sigma_2$), we construct the net with base point
$W^{\dS}_{x_2}$. We consider the net $\sN^{\dS}$ of standard subspaces on de Sitter spacetime
defined as in \eqref{eq:exgen1} by  
\begin{equation}
  \label{eq:hds4}
  \cW_+\ni W^{\dS}\mapsto {\sN^{\dS}(W^{\dS})}\subset\cH_U
  \quad \mbox{ defined by } \quad
  \sN^{\dS}(\Lambda(g)W^{\dS}_{x_2}):=V(g)\sN_U(W_{h}).   
\end{equation}
By Theorem~\ref{thm:3.5}, this net  is Lorentz
  covariant but not modular covariant.
  To match notation, take $h_1:=k_1$ and $h_2:=h$.

  \begin{remark}
Consider the previous identification of $\SL_2(\RR)\cong H^\up 
\subset\PSL_n(\RR)$,
where $h = h_c -k_1$, and $k_1$ and $k_2$ 
are the Euler elements in $\fh \cong \fsl_2(\R)$
associated to the wedge domains
$W_{x_2}^{\dS}$ and $W_{x_1}^{\dS}$ in~$\dS^2$, respectively. 
Then 
\[ k_1 =r(\pi/2) k_2
r(\pi/2)^{-1}, \quad \mbox{  where } \quad 
r(\theta)= \pmat{ \cos(\theta/2) & \sin(\theta/2) \\ 
- \sin(\theta/2) & \cos(\theta/2)}.\]
 By the BGL construction and \eqref{eq:R}, 
\[ U(e^{-2\pi it \cdot h})=U(e^{-2\pi it \cdot (h_c - k_1)})=\Delta_{\sN^{\dS}(W_{x_2}^{\dS})}^{it} 
\quad \mbox{ and } \quad 
U(e^{-2\pi it \cdot (h_c - k_2)})=\Delta_{\sN^{\dS}(W_{x_1}^{\dS})}^{it}.\]
As the two Euler elements $h = h_c - k_1$ and $h_c - k_2$ generate
$\fb \cong \gl_2(\R)$, we also observe that 
modular covariance manifestly fails
because these two one-parameter groups
do not generate a representation of $\widetilde{\SL_2}(\R)$ (cf.~\cite[Thm.~1.7]{BGL95}).
\end{remark}

\subsubsection{{Counterexamples on Minkowski spacetime}}\label{sect:minkcount}

Let $\sN$ be a net of standard subspaces:
$$ {\cW^{\dS}_+} \ni W^{\dS}\longmapsto \sN(W^{\dS})\subset\cH$$
{which is covariant for a unitary} 
representation $V$ of the double covering 
$\SL_2(\RR)\simeq\tilde\cL_+^\up$  of the 
Lorentz group on the Hilbert space~$\mathcal H$. We assume that the couple $(V, \sN)$ satisfies \eqref{eq:Lorcov}.
We further assume that the net $\sN$ does not satisfy modular covariance, {i.e., that there exist
$g \in \SL_2(\R)$ and $t \in \R$ such that
\begin{equation}
  \label{eq:nomodcov}
  \Delta_{\sN(W^{\dS}_{x_1})}^{-it} \sN(g W^{\dS}_{x_1})
  \not
  = \sN(\Lambda_{W_{x_1}^{\dS}}(2\pi t)g W^{\dS}_{x_1}) 
  = \sN(\Lambda(\exp(th_1))g W^{\dS}_{x_1}).  
\end{equation}
(Note that the two wedges $W^{\dS}_{x_1}$ and $W^{\dS}_{x_2}$ are Lorentz conjugate, so
that the base point in the wedge space does not matter.)
We constructed such nets in Sections {\ref{sect:desittcex} and} \ref{ex:SLn}.

The proper Poincar\'e group
\[ \cP_+=\cP_+^\up\cup\cP_+^\down \quad \mbox{ with }\quad 
\cP_+^{\up\down}=\RR^{1+2}\rtimes \cL_+^{\up\down} \]
is the inhomogeneous Lorentz group acting on the Minkowski spacetime:
the real space $\RR^{1+2}$ endowed with the Lorentzian metric with signature $(1,-1,-1)$. We write
\[ \cW^{\RR^{1+2}}_+=\cP_+^\up.W_{x_1} \]
for the set of wedges of $(1+2)$-dimensional Minkowski spacetime.
Note that ${\dS}^2\subset \RR^{1+2}$ and that, 
for every $W^{\dS}\subeq \dS^2$, there exists
a unique $W^{\RR^{1+2}}\in \cL_+^\up.W_{x_1}^{\R^{1+2}} \subeq \cW_+^{\RR^{1+2}}$
such that $W^{\dS}=W^{\RR^{1+2}}\cap \dS^2$ (cf.~Remark~\ref{ex:desit}).
We may thus identify $\cW_+^{\dS}$ with a
subset of $\cW_+^{\RR^{1+2}}$. The Poincar\'e group
$\cP_+^\uparrow=\RR^{1+2}\rtimes\cL_+^\up$ act on $\cW_+^{\dS}$
through the projection onto $\cL_+^\up$, so that
the translation group acts trivially.

With this identification, $\sN$ extends to a net $\hat\sN$
on Minkowski wedges by 
\[ \cW_+^{\RR^{1+2}}\ni W^{\RR^{1+2}}\mapsto
\hat\sN(W^{\RR^{1+2}})\subset\mathcal H \] 
with
\[ \hat\sN\big(\Lambda^{\RR^{1+2}}(a,g)W_{x_1}^{\RR^{1+2}}\big)
:=V(g)\sN(W_{x_1}^{\dS}),\qquad (a,g)\in\tilde\cP_+^\up\]
and $\Lambda^{\RR^{1+2}}:\tilde \cP_+^\uparrow=\RR^{1+2}\rtimes \tilde\cL_+^\uparrow\ni (a,g)\longmapsto(a,\Lambda(g))\in\cP_+^\uparrow$  is the quotient map.
As above, this construction requires that the stabilizer
subgroup of the standard right wedge
$W_R := W^{\R^{1+2}}_{x_1}$ in $\cP_+^\up$ preserves
$W^{\dS}_{x_1}$ (\cite[Lemma~4.13]{NO17}).

Now let $U_0$ be an (anti-)unitary positive energy\footnote{Positive energy means that the joint spectrum of the translation generator lies inside the forward light cone $C=\{x\in\RR^{1+2}:x^2\geq0, x_0\geq0\}$} representation of
$\cP_+$ on the Hilbert space $\mathcal  K$.
We write $\sN_{U_0}$ for the
BGL-net of standard subspaces associated to $U_0$ (Theorem \ref{thm:MN}).
We define a unitary positive energy  representation of   $\tilde \cP_+^\up=\RR^{1+2}\rtimes\SL_2(\RR)$, the double covering of $\cP_+^\up$, on the tensor product Hilbert space $\mathcal H\otimes \mathcal K$ by 
\begin{equation}\label{eq:UUU}
 U: \tilde{\mathcal P}_+^\up=\RR^{1+2}\rtimes \SL_2(\RR)\ni 
(a,g)\longmapsto V(g)\otimes U_0(a,g)\in \U(\mathcal H\otimes \mathcal K). \end{equation}
The tensor product subspaces  $\hat\sN(W^{\RR^{1+2}})\otimes \sN_{U_0}(W^{\RR^{1+2}})\subset\mathcal H\otimes\mathcal K$  
are standard by  \cite[Prop.~2.6]{LMR16}.
We obtain a (non-local) net of standard subspaces on Minkowski spacetime
{
\[ \sH:\;\cW^{\RR^{1+2}}_+\ni W^{\RR^{1+2}}\longmapsto 
\sH(W^{\RR^{1+2}})\subset \mathcal H\otimes \mathcal K,\]}
where $\sH(W_R)=\hat\sN(W_R)\otimes \sN_{U_0}(W_R)$ and, for the standard right wedge
$W_R := W_{x_1}^{\R^{1+2}}$: 
\begin{align*}\sH(W^{\RR^{1+2}}) &=\sH(\Lambda^{\RR^{1+2}}((a,g)).W_R)
 =U((a,g))\sH(W_R)\\&=(V(g)\otimes U_0(a,g))(\hat\sN(W_R)
  \otimes \sN_{U_0}(W_R))\quad \mbox{ for }
  \quad (a,g)\in \tilde\cP_+^\up = \RR^{1+2}\rtimes \SL_2(\RR).
\end{align*}

Poincar\'e covariance is satisfied by construction, 
but modular covariance is broken.
Indeed, by  \cite[Sect.~2]{LMR16}, $$\Delta^{it}_ {\hat\sN(W)\otimes \sN_{U_0}(W)}=\Delta^{it}_{\hat\sN(W)}\otimes \Delta^{it}_{\sN_{U_0}(W)}$$ and since $\hat\sN$ does not satisfy 
modular covariance, there exists a $g\in\tilde\cL_+^\up \subeq \tilde{\cP}_+^\up$ such that
we obtain: 
\begin{align*}
\Delta^{-it}_{\sH(W_R)}\sH(gW_R)
  &=(\Delta^{-it}_{\hat\sN(W_R)}\otimes \Delta^{-it}_{\sN_{U_0}(W_R)})(\hat\sN(gW_R)\otimes
    \sN_{U_0}(gW_R))\\
&=(\Delta^{-it}_{\hat\sN(W_R)}\hat\sN(gW_R))\otimes \Delta^{-it}_{\sN_{U_0}(W_R)}\sN_{U_0}(gW_R)\\
&=(\Delta^{-it}_{\hat\sN(W_R)}\hat\sN(gW_R))\otimes U_0(\Lambda_{W_R}(2\pi t))\sN_{U_0}(gW_R)\\
&=(\Delta^{-it}_{\hat\sN(W_R)}\hat\sN(gW_R))\otimes \sN_{U_0}(\Lambda_{W_R}(2\pi t)gW_R)\\
&\not=\hat\sN(\Lambda_{W_R}(2\pi t)gW_R)\otimes\sN_{U_0}(\Lambda_{W_R}(2\pi t)gW_R)
\end{align*}
where the last inequality follows from \eqref{eq:nomodcov}.
}

\section{Outlook: Disintegration, locality and higher dimensions}
\mlabel{sec:4} 

{\bf Comments on the representation theory.} Let $\sN_U$
be the BGL-net associated to an anti-unitary (positive energy) representation $U$ of a $\mathbb Z_2$-graded Lie group $G$ whose Lie algebra $\fg$ contains Euler elements, cf.~\cite{MN21}. The counterexamples to the BW property described in \cite{Mo18, LMR16} (see also \cite[Thm.~3.1]{Bo98})  can be interpreted in our general setting  in the following sense:
 
 \begin{proposition} \mlabel{lem:4.6} 
Let $(U,\cH)$ be an (anti-)unitary representation 
of the graded group $G$ and 
\[ \zeta \: G \to U(G)' \cap \U(\cH) \] 
be a group homomorphism. Then 
\[ \tilde U(g) := U(g) \zeta(g) = \zeta(g) U(g)  \] 
defines an (anti-)unitary representation of $G$ on $\cH$, and if 
$\sN := \sN_U \: \cG(G)\to \Stand(\cH)$ is 
the BGL net {corresponding to $U$}, 
then $\sN$ is $\tilde U$-covariant and satisfies the 
modular covariance condition.
\end{proposition}

\begin{prf} The net $\sN$ is modular covariant by construction and
    since each $\zeta(g)$ fixes all subspaces in the net by Lemma~\ref{lem:sym},
    the $\sN$ is $\tilde U$-covariant. 
\end{prf}
In the setting of Proposition 4.1, if $U$ is an infinite multiple of an irreducible Poincar\'e representation, then $\tilde U$ disintegrates with infinitely many disjoint Poincar\'e representations, see e.g.~\cite[Sect.~5]{Mo18}, \cite[Sect.~7]{LMR16}.

In the current paper,  our family of counterexamples to the modular covariance property on de Sitter spacetime relies on the inclusion
\[ H^\up=\SL_2(\RR)\subset G:=\GL_2(\RR)\qquad  \text{with}\qquad \fh:=\fsl_2(\RR)\subset\fg=\fgl_2(\RR),\]
where we consider a {non-symmetric Euler element
  $h\in \fgl_2(\RR)$ and observe that
  $\fg^h\not\supseteq[\fg,\fg]=\fsl_2(\RR)$}.
Note that $\fgl_2(\RR)=\fsl_2(\RR)\oplus \fz(\fg)$ where $\fz(\fg)=\RR \cdot\textbf{1}$.

Let $(V,\sH)$ be a net of standard subspaces  satisfying \eqref{eq:Lorcov} on the  two-dimensional de Sitter spacetime,
where $V$ is a  representation of the (covering of the) Lorentz group acting covariantly on $\sH$. Assume that $\sH$ has been
constructed as in Section~\ref{sect:desittcex}, 
with respect the group $G=\GL_2(\RR)$.
Then $h$ is a  non-symmetric Euler element in $\fg=\gl_2(\RR)$
as above and {$H^\up = \la \exp \fh \ra =\SL_2(\RR)\subset\GL_2(\RR)$}.

An irreducible  unitary representation of $G^\up=\R^\times_+\SL_2(\RR)$ is of the form $U_p(g,a)= V(g)e^{ipa}$ where $V$ is an irreducible representation of $\SL_2(\RR)$ and $\RR\ni p\neq0$. The representation 
 $V$ extends to an (anti-)unitary irreducible representation of
 $H=H^\up\rtimes\mathbb Z_2$ on the same space (see e.g.~\cite[Thm.~2.24]{MN21}).
  An (anti)-unitary representation of $G=G^\up\rtimes\mathbb Z_2$
  which is non-trivial on the center $\textbf{1}\times \RR_+^\times\subset G^\up$
  decomposes on $G^\up$ as  $U:=U_p\oplus U_{-p}$ because of the antilinearity of the action of the involution.

Let $V=U|_H$ and $\sH(W_{x_1}^{\dS})=\sN_U(W_h)$. So the net $\sH(gW_{x_1}^{\dS})=V(g)\sH(W_h)$ is defined by covariance but it does not satisfy modular covariance. We conclude that, on de Sitter spacetime, even if $V$ is a two-fold direct sum of irreducible representations of $\SL_2(\RR)$, it is possible that a  Lorentz covariant net of standard subspaces $(V,\sH)$ is not modular covariant. On the other hand, on Minkowski spacetime, due to the tensor product representation \eqref{eq:UUU}, the Poincar\'e representation $U$ contains infinitely many inequivalent representations
in the direct integral {decomposition} (cf.~\cite[Sect.~5]{Mo18}, \cite[Sect.~7]{LMR16}).\\

\nin{\bf Counterexamples with bosonic representations.} Lemma \ref{lem:slgl} claims that if $\fg$ is simple then $\Inn_\fg([\fb,\fb])\simeq \SL_2(\RR)$. Then, going back to the construction presented in Section \ref{sect:desittcex}, the covariant representation of the symmetry group acting on the two-dimensional de Sitter spacetime is actually a representation of the double covering of the Lorentz group that does not factor through the Lorentz group. In order to consider representations
  of the Lorentz group, one can get rid of the assumption on $\fg$ to be simple. Note that Lemma \ref{lem:slgl}(a)-(c) do not require $\g$ to be simple.
   
  With the notation of Section \ref{sect:desittcex}, consider a (non simple) Lie algebra $\fg =\fsl_2(\RR)\oplus \R \xi$,  where $ \fz(\g)=\RR\cdot\xi$ is the one-dimensional center of $\fg$. Then $h_2=\xi-h_1\in \fg$  is a non-symmetric Euler element of $\fg$ and $h_1$ is an Euler element in $\fh=\fsl_2(\RR)$.
  It is clear that $[\fh,h_1-h_2]\not\subseteq\ker (h_2)$. 
  Let $G^\up=\PSL_2(\RR)\times \RR$ and $H^\up=\PSL_2(\RR)\times \{0\}$,
  then $\mathrm{Lie}(G^\uparrow)=\fg$ and $\mathrm{Lie}(H^\up)=\fh$. Consider the $\mathbb Z_2$-graded extension $H$ of $H^\up$ obtained adding an Euler involution.
  Then $G=H\times \RR$ is a $\mathbb Z_2$-graded extension of $G^\up$.  Let $U$ be an (anti-)unitary representation of $G$ on a Hilbert space $\cH_U$ and $\sN_U$ be the BGL-net with respect to the wedge set $\cW_+(W_{h_2})$. One can replicate the construction of the non-modular covariant net presented in
  Section~\ref{sect:desittcex} since $H_{h_1}=\exp(\RR h_1)$, and obtain a net of standard subspaces indexed by wedge regions on two-dimensional de Sitter spacetime $$\cW^{\dS}\supset W^{\dS}\longmapsto\sN(W^{\dS})\subset\cH_U.$$ The net $\sN$ is Lorentz covariant with respect to $V=U|_H$ as in $\eqref{eq:Lorcov}$,
  where $V$ {is already defined on $\mathcal L_+^\up$ and not only on   the $2$-fold covering.} \\

  \nin \textbf{Twisted-local non modular Lorentz covarariant nets on de Sitter spacetime. }The counterexamples  to the modular covariance property we provided in this paper are not compatible with locality {since the Euler element we start the construction with is not symmetric}. Indeed, the subspaces corresponding to causal complementary wedges  (cf.\ Definition~\ref{def:euler}(d))  are not local
  in the sense that, if $W_1\subset W_2'$, then  $\sH(W_1)\subset\sH(W_2)'$.
  Clearly, on de Sitter spacetime, since there are no wedge inclusions,
  the locality codition becomes  $\sH({W^{\dS}}')\subset\sH(W^{\dS})'$.
On the other hand we can construct twisted-local nets that are not modular covariant as follows.
 
We refer to the notation contained in Section~\ref{sect:desittcex}. 
Let $(\sN,U)$ be a non-modular covariant net  on de Sitter spacetime   as constructed in Section~\ref{sect:desittcex},  with {$\sN(W_{x_1})=\sN_U(W_{h})$} where $h$ is  a non-symmetric Euler element in $\fg$.
We can define a second Lorentz covariant net of standard subspaces
on de Sitter spacetime  by covariance  under $U$,
starting  from
$$ \sK(W_{x_1}'^{\dS}):= \sN_U(W_{-h})=\sN_U(W_h)'\subset\cH,$$  and putting  
\[ \sK(
  \Lambda(g){W_{x_1}^{\dS}}'):=U(g)\sN_U(W_{h})'. \] 
Note that ${W_{x_1}^{\dS}}'=R(\pi)W_{x_1}^{\dS}$,
so  $\sK(W_{x_1}^{\dS})=U(r(\pi))\sN({W_{x_1}^{\dS}})'$, and by covariance \begin{equation}\label{eq:prime}
\sN({W^{\dS}}')=\sK(W^{\dS})'.
\end{equation}
The net is Lorentz covariant \eqref{eq:Lorcov} since, 
for $g\in H^\up=\SL_2(\RR)$, $\Lambda(g)W^{\dS}=W^{\dS}$ is equivalent to
$\Lambda(g){W^{\dS}}'={W^{\dS}}'$. 
With the notation of Section~\ref{sect:constr1},
{we put $h_2:=h$ and $h_1:=h_{W_{x_1}}$, so that \eqref{eq:negcond} holds:}
$$[\fh,h_1-h_2]\not\subseteq \ker(h_2).$$ 
In particular, for
$h_{W_h'}=-h_2$  
and $h_{W_{x_1}'}=-h_1$, we get  again $[\fh,h_2-h_1]\not\subseteq \ker(h_2),$ hence  $\sK$ is not modular covariant.
 The net $$\cW^{\dS}_+\ni W^{\dS}\longmapsto\tilde \sH(W^{\dS})= \sN(W^{\dS})\oplus \sK(W^{\dS})\subset \cH\oplus\cH$$ is Lorentz covariant with respect to the representation $$\SL_2(\RR)\ni g\longmapsto\tilde U(g)=U(g)\oplus U(g).$$
 Consider the operator $$Z=\left(\begin{array}{cc}\textbf{0}&\textbf{1}\\\textbf{1}&\textbf{0}
                                 \end{array}\right) \in \tilde U(\SL_2(\RR))'.$$
 The net $\tilde \sH$ is \textbf{twisted Haag dual} (cf.~{[GL95, MN21]}) with respect to the twist operator $Z$, in the sense that
\begin{equation}\label{eq:twisthaag}
  \tilde \sH(W')=Z\tilde\sH(W)'.
\end{equation}
 Indeed, {\eqref{eq:twisthaag} follows from} 
\begin{align*}
\tilde H({W_{x_1}^{\dS}}')=\tilde U(r(\pi))\tilde\sH({W_{x_1}^{\dS}})&=U(r(\pi))\sN(W_{x_1}^{\dS})\oplus U(r(\pi))\sK(W_{x_1}^{\dS})\\
&=U(r(\pi))\sN(W_{x_1}^{\dS})\oplus U(r(\pi))U(r(\pi))\sN(W_{x_1}^{\dS})'\\
&=\sN({W_{x_1}^{\dS}}')\oplus U(r(2\pi))\sN(W_{x_1}^{\dS})'\\
&\stackrel{\eqref{eq:prime}}=\sK(W_{x_1}^{\dS})'\oplus\sN(W_{x_1}^{\dS})'\\
&=Z(\sN(W_{x_1}^{\dS})'\oplus\sK(W_{x_1}^{\dS})')\\
&=Z\tilde\sH({W_{x_1}^{\dS}})'.
\end{align*}
By covariance  we obtain
$\tilde \sH({W^{\dS}}')=Z\tilde\sH(W^{\dS})'$ for every wedge on de Sitter spacetime.
 An analogous example can be constructed on Minkowski spacetime.\\

 \nin {\bf Higher dimensional spacetimes.} In this paper we {provide} counterexamples to modular covariance on two-dimensional
 de Sitter spacetime and three-dimensional Minkowski spacetime. The general construction principle
 (Section~\ref{sect:constr1}) {depends} neither on the spacetime dimension nor on the manifold we consider. Here we sketch how to extend the concrete construction presented in Section~\ref{ex:SLn}
 to higher dimensional de Sitter or Minkowski spacetimes. Let $\fh:=\so_{1,n}(\R)
 \subset \fg:=\fsl_m(\RR)$ and $h$ be a non-symmetric Euler element of $\fsl_m(\RR)$ such that
 {$\{0\} \not=  [\so_{1,n}(\R),h]\subeq \so_{1,n}(\R)$}
 and let $U$ be an (anti-)unitary representation of the $\Z_2$-graded
 extension $\PGL_m(\R)$ of $\PSL_m(\R)$. 
 Let $W_{x_1}^{\dS,\RR^{1+n}}$ be the wedge in the $x_1$-direction on de Sitter or on Minkowski spacetime, with
 the identification $\sH(W_{x_1})={\sN_U(W_h)}$. 
 One can define by covariance a Lorentz \eqref{eq:Lorcov} or a  Poincar\'e covariant \eqref{eq:Poicov} net of standard subspaces  if and only if
 \[ \Stab_{	\widetilde{\mathcal{L}}_+^\up} \,W_{x_1}=\{g\in\widetilde\cL_+^\up: \Lambda(g)W_{x_1}=W_{x_1}\}
   \subset (\widetilde\cL_+^\up)^h\]
or
 \[ \Stab_{\widetilde{\mathcal{P}}_+^\up} \,W_{x_1}=\{g\in\widetilde\cP_+^\up: \Lambda(g)W_{x_1}=W_{x_1}\}\subset (\widetilde{\cP}_+^\up)^h,\]
 respectively. 
In this case, the covariance conditions \eqref{eq:Lorcov} and \eqref{eq:Poicov} defines {standard}  subspaces
 \[ \sH(\Lambda(g)W_{x_1}^{\dS,\RR^{1+s}}):=\sN_U(gW_h)\quad \mbox{ for } \quad  g\in\widetilde{\cL}_+^\up \quad\text{or} \quad g\in\tilde \cP_+^\up,\]
 where $\widetilde{\mathcal L}_+^\up$ and $\widetilde{\mathcal P}_+^\up$ are the double (and universal) covering of the Lorentz and the Poincar\'e group respectively,  and $\Lambda$ is the covering homomorphism.
 
 \bigskip 
 
\nin{ \bf Acknowledgements.} VM was supported by Alexander-von-Humboldt Foundation through a Humboldt Research Fellowship for Experienced Researchers. VM also thanks the Mathematics Department of Sapienza University of Rome and acknowledges the University of Rome ``Tor Vergata'' funding OAQM, CUP E83C22001800005 and the European Research Council Advanced Grant 669240 QUEST.

KHN acknowledges support by DFG-grant NE 413/10-1.

\end{document}